\documentclass[preprint,12pt]{article}
\UseRawInputEncoding
\pdfoutput=1
\usepackage{graphicx}
\usepackage{float}
\usepackage{amssymb}
\usepackage{mathrsfs}
\usepackage{enumerate}
\usepackage{colortbl,dcolumn}
\usepackage{tabularx}
\usepackage{amsmath}
\usepackage{psfrag}
\usepackage{booktabs}
\usepackage{amsthm}
\usepackage{bm}
\usepackage{framed,multirow}
\usepackage{latexsym}
\usepackage{url}
\usepackage{lineno,hyperref}
\usepackage{epsfig}
\usepackage{multirow}
\usepackage{amsthm}
\usepackage{hyperref}
\usepackage{cleveref}
\usepackage{caption}
\usepackage{amsfonts}
\usepackage{float}
\usepackage{subcaption} 
\usepackage{etoolbox}
\usepackage[utf8]{inputenc}
\usepackage[T1]{fontenc}
\usepackage[comma, numbers]{natbib}
\usepackage[top=1in, bottom=1in, left=0.7in, right=0.7in]{geometry}

\newtheorem{remark}{Remark}
\newtheorem{assumption}{Assumption}
\newtheorem{theorem}{Theorem}
\newtheorem{lemma}{Lemma}
\newtheorem{definition}{Definition}
\newtheorem{proposition}{Proposition}
\begin{document}
\title{Parareal algorithms for stochastic Maxwell equations with the damping term driven by additive noise}
       \author{
      Liying Zhang\footnotemark[1],
        Qi Zhang\footnotemark[2]\\
       {\small
       \footnotemark[1]~\footnotemark[2]~School of Mathematical Science, China University of Mining and Technology, Beijing 100083, China.}
    }
       \footnotetext{\footnotemark[2]Corresponding author: zq17866703992@163.com}
\maketitle
\begin{abstract}    
 In this paper, we investigate the strong convergence analysis of parareal algorithms for stochastic Maxwell equations with the damping term driven by additive noise. The proposed parareal algorithms proceed as two-level temporal parallelizable integrators with the stochastic exponential  integrator as the coarse  $\mathcal{G}$-propagator and both the exact solution integrator and the stochastic exponential integrator as the fine  $\mathcal{F}$-propagator. It is proved that the convergence order of the proposed algorithms linearly depends on the iteration number $k$. Numerical experiments are performed to illustrate the convergence of the parareal algorithms for different choices of the iteration number $k$ and the damping coefficient $\sigma$.\\

\textbf{Key Words: }{\rm\small} Stochastic Maxwell equations, Parareal algorithm, Strong convergence, Stochastic exponential integrator
\end{abstract}
\section{Introduction}
 When the electric and magnetic fluxes are perturbed by noise, the uncertainty and stochasticity can have a subtle but profound influence on the evolution of complex dynamical systems \citep{Rytovetal1989}. In order to model the thermal motion of electrically charged microparticles, we consider the stochastic Maxwell equations with damping term driven by additive noise as follows
\begin{align}	\label{eq1}
	\left\{
	\begin{array}{ll}
		\varepsilon \partial_t \bm{E}(t,\bm{x})=\nabla \times \bm{H}(t,\bm{x})-\sigma\bm{E}(t,\bm{x})-J_e (t,x,\bm{E},\bm{H}) -J_{e}^r(t,\bm{x}) \cdot \dot{W}	, & (t,\bm{x})\in \left(0,T \right] \times D,\\
		\mu \partial_t \bm{H}(t,\bm{x})=-\nabla \times \bm{E}(t,\bm{x})-\sigma\bm{H}(t,\bm{x})-J_m (t,x,\bm{E},\bm{H}) -J_{m}^r(t,\bm{x}) \cdot \dot{W}	,  &(t,\bm{x})\in \left(0,T \right] \times D, \\
		\bm{E}(0,\bm{x})=\bm{E}_{0}(\bm{x}),
		\bm{H}(0,\bm{x})=\bm{H}_{0}(\bm{x}), &\bm{x}\in D, \\
		\bm{n} \times \bm{E}=0, &(t,\bm{x})\in \left(0,T \right] \times \partial D,
	\end{array}
	\right.
\end{align}
where $D\subset\mathbb{R}^{3}$ is an open, bounded and Lipschitz domain with boundary $\partial D$, of which $\bm{n}$ is the unit outward. Here $\varepsilon$ is the electric permittivity and $\mu$ is the magnetic permeability. The damping terms $\sigma\bm{E}$ and $\sigma\bm{H}$ are usually added to simulate the attenuation of electromagnetic waves in the medium, which can be caused by absorption, scattering or other non-ideal factors in the medium. The function $J_e$ and $J_{e}^r$ describe electric currents (or $J_m$ and $J_{m}^r$ describe magnetic currents). In particular, $J_{e}^r$ and $J_{m}^r$ do not depend on the electromagnetic fields $\bm{E}$ and $\bm{H}$. The authors  in \citep{Liaskosetal2010} proved the mild, strong and classical well-posedness for the Cauchy problem of stochastic Maxwell equations.  Meanwhile, the authors in \citep{Horsinetal2010} studied the approximate controllability of the stochastic Maxwell equations via an abstract approach and a constructive approach using a generalization of the Hilbert uniqueness method. Subsequently the work \citep{Roachetal2012} combined the study of well-posedness, homogenization and controllability of Maxwell equations with the description of the constitutive relations of complex media and dealt with deterministic and stochastic issues in both the frequency and time domains. 

Since stochastic Maxwell equations are a kind of stochastic Hamiltonian PDEs, constructing stochastic multi-symplectic numerical  methods for problem (\ref{eq1}) has been paid  more and more attention. The stochastic multi-symplectic numerical method for stochastic Maxwell equations driven by additive noise was proposed in \citep{Hongetal2014} based on the stochastic variational principle. Subsequently the authors in \citep{Chenetal2016} used a straightforward approach to avoid the introduction of additional variables and  obtained three effecitve stochastic multi-symplectic numerical methods. Then the authors in \citep{Hongetal2017}  used the wavelet collocation method in space and the stochastic symplectic method in time to	construct the stochastic multi-symplectic energy-conserving method for three-dimensional stochastic Maxwell equations driven by multiplicative noise. The work in \citep{Zhangetal2019a} made a review on these stochastic multi-symplectic methods and summarised numerical methods for various stochastic Maxwell equations driven by additive and multiplicative noise.  The general case of stochastic Hamiltonian PDEs was considered in \citep{ZhangJi2019}, where the multi-symplecticity of stochastic RK methods was investigated. Recently, the authors in \citep{Sunetal2022} and \citep{Sunetal2023} constructed
multi-symplectic  DG methods for stochastic Maxwell equations driven by additive noise and multiplicative noise.
Furthermore, the work in \citep{Hongetal2022} employed the local radial basis function collocation method and the work in \citep{Hou2023} utilized the global radial basis function collocation method for stochastic Maxwell equations driven by multiplicative noise to preserve  multi-symplectic structure. Additionally, \citep{Chen2021} developed a symplectic discontinuous Galerkin full discretisation method for stochastic Maxwell equations driven by additive noise. Other efficient numerical methods for stochastic Maxwell equations also are investigated, see \citep{Zhang2008} for the finite element method, \citep{Badieirostamietal2010} for the numerical method based on the Wiener chaos expansion, \citep{Chenetal2022} for ergodic numerical method, \citep{Chenetal2021} for operator splitting method and  \citep{ZhouLiang2024} for CN-FDTD and Yee-FDTD methods. 

Meanwhile, there are a lot of pregnant works focused mainly on strong convergence analysis of the numerical methods for stochastic Maxwell equations. In the temporal discretization methods, the semi-implicit Euler method was proposed in \citep{Chenetal2019b} to proved mean-square convergence order is $1/2$ for stochastic Maxwell equations driven by multiplicative noise. Subsequently the work in \citep{Chenetal2019a} studied the stochastic Runge-Kutta method with mean-square convergence order 1 for stochastic Maxwell equations driven by additive noise. In addition, explicit exponential integrator was proposed in \citep{Cohenetal2020} for stochastic Maxwell equations   with mean-square convergence order $1/2$ for  multiplicative noise and convergence order $1$ for additive noise. The work \citep{Chen2021} developed discontinuous Galerkin full discretization method for stochastic Maxwell equations driven by additive noise with mean-square convergence order $k/2$ in time and $k-1/2$ in space, where $k=1,2$ represents $H_{k}$ regularity.	Another related work by authors of \citep{Chenetal2022} showed the ergodic discontinuous Galerkin full discretization for
stochastic Maxwell equations with mean-square convergence order both $1/2$ in
the temporal and spatial directions. In recent works
\citep{Sunetal2022} and \citep{Sunetal2023}, high order discontinuous Galerkin methods were designed for the stochastic Maxwell equations driven by  additive noise and  multiplicative noise with mean-square convergence order both $k+1$. Besides, the authors of \citep{Chenetal2021} presented the operator splitting method for stochastic Maxwell equations driven by additive noise with mean-square convergence order $1$.

In order to increase the convergence order and improve the computational efficiency on stochastic differential equations, the parareal algorithm has received attentions. This algorithm we focus on is a two-stage time-parallel integrator originally proposed in \citep{Lionsetal2001} and further works studied on theoretical analysis and applications for differential model problems, see, for instance, \citep{BalMaday2002,WuZhou2015,DaiMaday2013,GanderVandewalle2007,GanderHairer2014,Daietal2013}. In terms of stochastic model, the work in
\citep{Zhangetal2019b} investigated  the parareal algorithm combining the projection method to SDEs with conserved quantities. Then the parareal algorithm for the stochastic Schr\"odinger equations with weak damping term driven by additive noise was studied in \citep{Hongetal2019} with fine propagator being the exact solver and coarse propagator being the exponential $\theta$-scheme. And the proposed algorithm increases the convergence order to $k$  in the linear case for $\theta \in[0,1]\backslash\{1/2\}$. The parareal algorithm for semilinear parabolic SPDEs behaved differently in \citep{BrehierWang2020} depending on the choice of the coarse integrator. When the linear implicit Euler scheme was selected, the convergence order was limited by the regularity of the noise with the increase of  iteration number, while for the stochastic exponential scheme, the convergence order always increased.
To the best of our knowledge, there has been no reference considering the convergence analysis of the parareal algorithm for stochastic Maxwell equations till now. 

Inspired by the pioneering works, we establish strong convergence analysis of the parareal algorithms for stochastic Maxwell equations with damping term driven by additive noise. Combining the benefits of the stochastic exponential integrator, we use this integrator as the coarse $\mathcal{G}$-propagator and for the fine  $\mathcal{F}$-propagator, two choices are considered: the exact solution integrator as well as the stochastic exponential integrator. Taking advantage of the contraction semigroup generated by the Maxwell operator and the damping term, we derive the uniform mean-square convergence analysis of the proposed parareal algorithms with convergence order $k$. The key point of convergence analysis is that  the
error between the solution computed by the parareal algorithm and the reference solution generated by the fine propagator for
the stochastic exponential integrator still maintains the consistent convergence results. Different from the exact solution integrator as the fine  $\mathcal{F}$-propagator, we need to make use of the Lipschitz continuity of the residual operator rather than the integrability of the exact solution directly in this case, which requires us to make assumptions about the directional derivatives of the drift coefficient. We find that the selection of parameters have an impact on the convergence analysis results of the parareal algorithms. An appropriate damping coefficient ensures stability and accelerates the convergence results and the scale of noise induces a perturbation of the solution numerically.

The  article is organized as follows. In the forthcoming section, we collect some preliminaries about stochastic Maxwell  equations. In section 3, we devote to introducing the parareal algorithms based on the exponential scheme as the coarse $\mathcal{G}$-propagator and both the exact solution integrator and the stochastic exponential integrator as the fine  $\mathcal{F}$-propagator. In section 4, two convergence results in the sense of mean-square are analyzed. In section 5, numerical experiments are dedicated to illustrate the convergence analysis with the influences on the iteration number $k$ and the damping coefficient $\sigma$ and the effect of noise with different scale $\lambda$ on the numerical solution.

 To lighten notations, throughout this paper,  C stands for a constant which might be dependent of $T$ but is independent of $\Delta T$ and may vary from line to line.
	\section{Preliminaries}
The basic Hilbert space is $\mathbb{H}:= L^{2}(D)^{3}\times L^{2}(D)^{3}$ with inner product 
\begin{align*}
	\left\langle
	\begin{pmatrix}
		\bm{E}_{1}\\
		\bm{H}_{1}
	\end{pmatrix},
	\begin{pmatrix}
		\bm{E}_{2}\\
		\bm{H}_{2}
	\end{pmatrix}
	\right\rangle
	_\mathbb{H}
	=\int_{D}(\varepsilon\bm{E}_{1} \cdot \bm{E}_{2}+\mu \bm{H}_{1}\cdot\bm{H}_{2})d\bm{x},
\end{align*}
for all $\bm{E}_{1},\bm{H}_{1},\bm{E}_{2},\bm{H}_{2}\in L^{2}(D)^{3} $ and  the norm 
\begin{align*}
	\left\|
	\begin{pmatrix}
		\bm{E}\\
		\bm{H}
	\end{pmatrix}
	\right\|
	_\mathbb{H}^{2}
	=
	\int_{D}(\varepsilon \|\bm{E}\|^{2}+\mu \|\bm{H}\|^{2})d\bm{x},\, \forall \bm{E},\,\bm{H}\in L^{2}(D)^{3}.
\end{align*}
In addition,  assume that $\varepsilon$ and $\mu$ are bounded and uniformly positive definite
functions: $\varepsilon,\,\mu \in L^{\infty}(D)$, $\forall \varepsilon,\,\mu >0$.

The $Q$-Wiener process $W$ is defined on a given probability space $(\Omega,\mathscr{F},P,\{\mathscr{F\}_{t}}_{t \in[0,T]})$ and can be expanded in a Fourier series 
\begin{align}\label{eq2}
	W(t)=\sum\limits_{n=1}\limits^{\infty}\lambda_{n}^{1/2}\beta_{n}(t)e_{n}, \,t\in[0,T],
\end{align}
where $\{\beta_{n}(t)\}_{n=1}^{\infty}$ is a sequence of independent standard real-valued Wiener processes and $\{e_{n}\}_{n=1}^{\infty}$ is a complete orthonormal system of $\mathbb{H}$ consisting of eigenfunctions of a symmetric, nonnegative and finite trace operator $Q$, i.e., $Tr\left(Q\right)<\infty$ and
$Qe_{n}=\lambda_{n}e_{n}$ with corresponding eigenvalue $\lambda_{n} \geq 0$.

The Maxwell operator is defined by 
\begin{equation}\label{eq3}
	M
	\begin{pmatrix}
		\bm{E}\\
		\bm{H}
	\end{pmatrix}
	:=
	\begin{pmatrix}
		0& \varepsilon^{-1}\nabla \times\\
		-\mu^{-1}\nabla \times&0
	\end{pmatrix}
	\begin{pmatrix}
		\bm{E}\\
		\bm{H}
	\end{pmatrix},
\end{equation}
with domain 
\begin{equation*}
	\mathcal{D}(M)=
	\begin{Bmatrix}
		\begin{pmatrix}
			\bm{E}\\
			\bm{H}
		\end{pmatrix}
		\in \mathbb{H}:M
		\begin{pmatrix}
			\bm{E}\\
			\bm{H}
		\end{pmatrix}
		=
		\begin{pmatrix}
			\varepsilon^{-1}\nabla \times \bm{H}\\
			-\mu^{-1}\nabla \times \bm{E}
		\end{pmatrix}
		\in \mathbb{H}, \,\bm{n} \times \bm{E} \bigg|_{\partial D}=0
		
	\end{Bmatrix}.
\end{equation*}
Based on the closedness of the operator $\nabla \times$, we have the following lemma.
\begin{lemma}\label{lemma1}\citep{Chenetal2023}
	The Maxwell operator defined in (\ref{eq3}) with domain $\mathcal{D}(M)$ is closed and skew-adjoint, and generates a  $C_{0}$-semigroup $S(t)=e^{tM}$on $\mathbb{H}$ for $t\in [0,T]$. Moreover, the frequently used property for Maxwell operator M is : $\langle Mu, u\rangle_{\mathbb{H}}=0$.
\end{lemma}
Let the drift term  $F:[0,T] \times \mathbb{H} \rightarrow \mathbb{H}$ be a Nemytskij operator associated with $\bm{J}_{e},\bm{J}_{m}$ defined by
\begin{equation*}
	F(t,u)(\bm{x})=
	\begin{pmatrix}
		-\varepsilon^{-1}\bm{J}_{e}(t,\bm{x},\bm{E}(t,\bm{x}),\bm{H}(t,\bm{x}))\\
		-\mu^{-1}\bm{J}_{m}(t,\bm{x},\bm{E}(t,\bm{x}),\bm{H}(t,\bm{x}))
	\end{pmatrix},
	\,\bm{x}\in D, \,u=(\bm{E}^{T},\bm{H}^{T})^{T} \in \mathbb{H}.
\end{equation*}
The diffusion term $B:[0,T] \rightarrow HS(U_{0}, \mathbb{H})$ is the Nemytskij operator defined by
\begin{equation*}
	(B(t)v)(\bm{x})=
	\begin{pmatrix}
		-\varepsilon^{-1}\bm{J}_{e}(t,\bm{x})v(\bm{x})\\
		-\mu^{-1}\bm{J}_{m}(t,\bm{x})v(\bm{x})
	\end{pmatrix},
	\,\bm{x}\in D, \,v \in U_{0}:=Q^{\frac{1}{2}}\mathbb{H}.
\end{equation*}
We consider the abstract form of (\ref{eq1}) in the infinite-dimensional space $\mathbb{H}:= L^{2}(D)^{3}\times L^{2}(D)^{3}$
\begin{align}
	\label{eq4}
	\left\{
	\begin{array}{ll}
		d u(t)=[Mu(t)-\sigma u(t)]dt+F(t,u(t))dt+B(t)dW,\,t\in\left(0,T \right],\\
		u(0)=u_{0},
	\end{array}
	\right.
\end{align}
where the solution $u=(\bm{E}^{T}, \bm{H}^{T})^{T}$ is  a stochastic process with values in $\mathbb{H}$. 

Let $\widehat S(t):=e^{ t(M-\sigma Id)}, \,\sigma \geq 0$ be the semigroup generated by operator $M-\sigma Id$. One can show that the damping stochastic Maxwell equations (\ref{eq4}) possess the following lemma.
\begin{lemma}\label{lemma2}
	For the semigroup $\{\widehat S(t)=e^{ t(M-\sigma Id)}, \,t\geq 0\}$ on $\mathbb{H}$, we obtain
	\begin{align*}
		\left\|\widehat S(t)\right\|_{\mathcal{L(\mathbb{H})}}\leq 1,
	\end{align*}	
	for all $t \geq 0$. 
\end{lemma}
\noindent{{\bf{Proof}.}} Based on the semigroup  $\{\widehat S(t)=e^{ t(M-\sigma Id)}, \,t\geq 0\}$ generated by the  operator $M-\sigma Id$, we deduce
\begin{align}\label{eq5}
	\left\|\widehat S(t)\right\|_{\mathcal{L(\mathbb{H})}}=\left\|e^{-\sigma t}S(t)\right\|_{\mathcal{L(\mathbb{H})}}\leq\left\| S(t)\right\|_{\mathcal{L(\mathbb{H})}}, \,t\geq 0.
\end{align}
Consider the deterministic system \citep{Chenetal2023}
\begin{align*}
	\left\{
	\begin{array}{ll}
		\dfrac{du(t)}{dt}=Mu(t), \,t\in\left(0,T \right],\\
		u(0)=u_{0}.
	\end{array}
	\right.
\end{align*}	
Thus
\begin{align*}
	\frac{d}{dt}\|u(t)\|_{\mathbb{H}}^{2}=2\left\langle\frac{du(t)}{dt}, u(t)\right\rangle_{\mathbb{H}}=2\left\langle Mu(t), u(t)\right\rangle _{\mathbb{H}}=0,
\end{align*}	
which leads to
\begin{align*}
	\|u(t)\|_{\mathbb{H}}=\|S(t)u_{0}\|_{\mathbb{H}}=\|u_{0}\|_{\mathbb{H}},
\end{align*}	
that is,
\begin{align*}
	\|S(t)\|_{\mathcal{L(\mathbb{H})}}=1.
\end{align*}	
Combining the formula (\ref{eq5}), we can conclude that the proof.\hfill $\Box$

To ensure the well-posedness of  mild solution of the stochastic Maxwell equations (\ref{eq4}), we need the following assumptions. 
\begin{assumption}\label{assump1}{\bf(Initial value)}.
	The initial value $u_{0}$ satisfies
	\begin{align*}
		\|u_{0}\|^{2}_{L_{2}(\Omega,\,\mathbb{H})}<\infty.
	\end{align*}
\end{assumption}
\begin{assumption}\label{assump2}{\bf(Drift nonlinearity)}.
	The drift operator $F$ satisfies
	\begin{align*}
		&\|F(t,u)\|_{\mathbb{H}} \leq C(1+\|u\|_{\mathbb{H}}),\\
		&\|F(t,u)-F(s,v)\|_{\mathbb{H}} \leq C(|t-s|+\|u-v\|_{\mathbb{H}}),
	\end{align*}
	for all $t,s \in [0,T], \,u,v \in \mathbb{H}$. Moreover, the nonlinear operator $F$ has bounded derivatives, i.e., 
	\begin{align*}
		\|DF(u).h\|_{\mathbb{H}} \leq C\|h\|_{\mathbb{H}},
	\end{align*}
	for $h\in \mathbb{H}$. 
\end{assumption}
\begin{assumption}\cite{Yan2005}\label{assump3}{\bf(Covariance operator)}.
	To guarantee the existence of a mild solution, we further assume the covariance operator $Q$ of $W(t)$ satisfies 
	\begin{align*}
		\|M^{(\beta-1)/2}Q^{1/2}\|_{\mathcal L_{2}(\mathbb{H})}<\infty, \quad\beta \in [0,1],
	\end{align*}
	where $\|\cdot\|_{\mathcal L_{2}(\mathbb{H})}$ denotes the Hilbert–Schmidt norm for operators from $\mathbb{H}$ to $\mathbb{H}$,\, $M^{(\beta-1)/2}$ is the $(\beta-1)/2$-th fractional powers of $M$ and  $\beta$ is a parameter
	characterizing the regularity of noise. In the article,
	we are mostly interested in $\beta=1$ for trace class operator $Q$. 
\end{assumption}
\begin{lemma}\label{lemma3}\citep{Chenetal2023}
	Let \Cref{assump1}, \Cref{assump2} and \Cref{assump3} hold, there exists a unique mild solution to ($\ref{eq4}$), which satisfies 
	\begin{equation*}
		u(t)=\widehat S(t)u_{0}+\int_{0}^{t}\widehat S(t-s)F(s,u(s))ds+\int_{0}^{t}\widehat S(t-s)B(s)dW(s),\,\mathbb{P} -a.s.,
	\end{equation*}
	for each $t\in[0,T]$, where $\widehat S(t)=e^{t(M-\sigma Id)}, \,t\geq 0$ is a $C_{0}$-semigroup generated by $M-\sigma Id$. 
	
	Moreover, there exists a constant $C\in (0,\infty)$ such that
	\begin{equation*}
		\sup\limits_{t\in[0,T]}\left\|u(t)\right\|_{L_{2}(\Omega,\, \mathbb{H})}\leq C(1+\left\|u_0\right\|_{L_{2}(\Omega,\, \mathbb{H})}).	
	\end{equation*}
\end{lemma}

The following lemma is the stability of analytical solution, which will be used in the proof of the Theorem  \ref{theorem1}.
\begin{lemma}\label{lemma4}\cite{Yan2005}
	If $u(t)$ and $v(t)$ are two solutions of ($\ref{eq4}$)  with different initial values $u_{0}$ and $v_{0}$, there exists a constant $C\in (0,\infty)$ such that
	\begin{equation*}
		\left\|u(t)-v(t)\right\|_{L_{2}(\Omega,\,\mathbb{H})} \leq C\left\|u_{0}-v_{0}\right\|_{L_{2}(\Omega,\,\mathbb{H})} .
	\end{equation*}
\end{lemma}

\section{Parareal algorithm for stochastic Maxwell equations}
\subsection{Parareal algorithm}
To perform the parareal algorithm, the considered interval $[0, T]$ is first divided into $N$ time intervals $[t_{n-1},t_{n}]$  with a uniform coarse step-size $\Delta T=t_{n}-t_{n-1}$ for any $n=1,\cdots,N$. Each subinterval is further divided into $J$ small time intervals $[t_{n-1,j-1},t_{n-1,j}]$ with a uniform fine step-size $\Delta t=t_{n-1,j}-t_{n-1,j-1}$ for all $n=1,\cdots,N$ and $j=1,\cdots,J$. The
parareal algorithm can be described as following
\begin{itemize}
	\item Initialization.	Use the coarse propagator $\mathcal{G}$ with the coarse step-size $\Delta T$ to compute initial value $u_{n}^{(0)}$ by
	\begin{align*}
		u_{n}^{(0)}&=\mathcal{G} (t_{n-1},t_{n},u_{n-1}^{(0)}),\quad n =1,\ldots ,N.\\
		u_{0}^{(0)}&=u_{0}.
	\end{align*}
	
	Let $K \in \mathbb{N}$ denote the number of parareal iterations: for all $n= 0,\ldots,N, \,k =0,\ldots,K-1$.
	\item Time-parallel computation.  Use the fine propagator $\mathcal{F}$ and  time step-size  $\Delta t$ to compute $\widehat u_{n}$ on each subinterval $[t_{n-1},t_{n}]$  independently
	\begin{align}\label{eq6}
		&\widehat u_{n-1,j}=\mathcal{F}(t_{n-1,j-1},t_{n-1,j},\widehat u_{n-1,j-1}), \quad j=1,\cdots,J ,\\
		&\widehat u_{n-1,0}=u_{n-1}^{(k)}.\nonumber
	\end{align}
	
	\item Prediction and correction. Note that we get two numerical solutions $u_{n}^{(0)}$ and $\widehat u_{n-1,J}$
	at time $t_{n}$ through the initialization and parallelization, the sequential
	prediction and correction is defined as
	\begin{align}\label{eq7}
		u_{n}^{(k+1)}&=\mathcal{G} (t_{n-1},t_{n},u_{n-1}^{(k+1)})+\widehat u_{n}-\mathcal{G} (t_{n-1},t_{n},u_{n-1}^{(k)}),\\
		u_{0}^{(k)}&=u_{0}.\nonumber
	\end{align}
	Noting  that equation (\ref{eq6}) is of the following form $\widehat u_{n}=\mathcal{F}(t_{n-1},t_{n},u_{n-1}^{(k)})$, then parareal algorithm can be written as 
	\begin{align}\label{eq8}
		u_{n}^{(k+1)}=\mathcal{G} (t_{n-1},t_{n},u_{n-1}^{(k+1)})+\mathcal{F}(t_{n-1},t_{n},u_{n-1}^{(k)})-\mathcal{G} (t_{n-1},t_{n},u_{n-1}^{(k)}).
	\end{align}
\end{itemize}
\begin{remark}
	The coarse integrator $\mathcal{G}$ is required to be easy to calculate and enjoys a less computational cost, but need not to
	be of high accuracy. On the other hand, the fine integrator $\mathcal{F}$ defined on each subinterval is assumed to be more accurate  but more costly than $\mathcal{G}$.  Note that $\mathcal{G}$ and $\mathcal{F}$ can be the same numerical method or different numerical methods. In the article, the exponential integrator is chosen as the coarse integrator $\mathcal{G}$ and both the exact integrator and the exponential integrator are chosen as the fine integrator $\mathcal{F}$.
\end{remark}
\subsection{Stochastic exponential scheme} 
Consider the mild solution of the stochastic Maxwell equations (\ref{eq4}) on the time interval $[t_{n-1},t_{n}]$ 
\begin{align}\label{eq9}
	u(t_{n})=\widehat S(\Delta T)u(t_{n-1})+\int_{t_{n-1}}^{t_{n}}\widehat S(t_{n}-s)F(u(s))ds+\int_{t_{n-1}}^{t_{n}}\widehat S(t_{n}-s)B(s)dW(s),
\end{align} 
where $C_{0}$-semigroup $\widehat S(\Delta T)=e^{\Delta T(M-\sigma Id) }$.

By approximating the integrals in above mild solution (\ref{eq9}) at the left endpoints, we can obtain the stochastic exponential scheme
\begin{align}\label{eq10}
	u_{n}=\widehat S(\Delta T)u(t_{n-1})+\widehat S(\Delta T)F(u(t_{n-1}))\Delta T+\widehat S(\Delta T)B(t_{n-1})\Delta W_{n},
\end{align}
where $\Delta W_{n}=W(t_{n})-W(t_{n-1})$.
\subsection{Coarse and fine propagators}
\begin{itemize}
	\item Coarse propagator.
	The stochastic exponential scheme is chosen as
	the coarse propagator  with
	time step-size $\Delta T$ by (\ref{eq10})
	\begin{align}\label{eq11}
		\mathcal{G}(t_{n-1},t_{n},u)=\widehat S(\Delta T)u+\widehat S(\Delta T)F(u)\Delta T+\widehat S(\Delta T)B(t_{n-1})\Delta W_{n},
	\end{align}
	where  $\widehat S(\Delta T)=e^{\Delta T(M-\sigma Id) }$ and $\Delta W_{n}=W(t_{n})-W(t_{n-1})$.
	\item Fine propagator.
	The exact solution as the fine propagator with time step-size $\Delta t$ by (\ref{eq9}) 
	\begin{align}\label{eq12}
		\mathcal{F}(t_{n-1,j-1},t_{n-1,j},u)=\widehat S(\Delta t)u+\int_{0}^{\Delta t}\widehat S(\Delta t-s)F(u(s))ds
		+\int_{0}^{\Delta t}\widehat S(\Delta t-s)B(t_{n-1,j-1})dW(s),
	\end{align}
	where  $\widehat S(\Delta t)=e^{\Delta t(M-\sigma Id)}$.
	
	Besides, the other choice is the stochastic exponential scheme is chosen as the fine propagator with time step-size $\Delta t$ by (\ref{eq10})
	\begin{align}\label{eq13}
		\mathcal{F}(t_{n-1,j-1},t_{n-1,j},u)=\widehat S(\Delta t)u+\widehat S(\Delta t)F(u)\Delta t+\widehat S(\Delta t)B(t_{n-1,j-1})\Delta W_{n-1,j},
	\end{align}
	where  $\widehat S(\Delta t)=e^{\Delta t(M-\sigma Id)}$ and $\Delta W_{n-1,j}=W(t_{n-1,j})-W(t_{n-1,j-1})$.
\end{itemize}
\section{Main results}
In this section, two convergence analysis results will be given, i.e., we investigate the parareal algorithms obtained by choosing the stochastic exponential integrator as the coarse integrator and both the exact integrator and the stochastic exponential integrator as the fine integrator.
\subsection{The exact integrator as the fine integrator $\mathcal{F}$}
\begin{theorem}\label{theorem1}
	Let \Cref{assump1}, \Cref{assump2} and \Cref{assump3} hold,  we apply the stochastic exponential integrator for coarse propagator $\mathcal{G}$ and exact solution integrator for fine propagator $\mathcal{F}$. Then we have the following convergence estimate for the fixed iteration number $k$
	\begin{equation}\label{eq14}
		\sup\limits_{1\leq n \leq N}\left\|u(t_{n})-u_{n}^{(k)}\right\|_{L_{2}(\Omega,\,\mathbb{H})} \leq \frac{C_{k}}{k!} \Delta T^{k} \prod \limits_{j=1}^k (N-j)\sup\limits_{1\leq n \leq N}\left\|u(t_{n})-u_{n}^{(0)}\right\|_{L_{2}(\Omega,\, \mathbb{H})},
	\end{equation}
	with a positive constant $C$ independent on $\Delta T$, where the parareal solution $u_{n}^{(k)}$ is defined in (\ref{eq8}) and the exact solution $u(t_{n})$ is defined in (\ref{eq9}).
\end{theorem}
To simplify the exposition, let us introduce the following notation.
\begin{definition}
	The residual operator 
	\begin{align}\label{eq15}
		\mathcal{R}(t_{n-1},t_{n},u):=\mathcal{F}(t_{n-1},t_{n},u)-\mathcal{G}(t_{n-1},t_{n},u),
	\end{align}
	for all $n\in{0,\cdots,N}$.
\end{definition}
Before the error analysis, the following two useful lemmas are introduced.
\begin{lemma}\label{lemma5}\cite{GanderVandewalle2007}
	Let $M:=M(\beta)_{N\times N}$ be a strict lower triangular Toeplitz matrix and its elements are defined as 
	\begin{align*}
		M_{i1}=
		\left\{
		\begin{array}{ll}
			0, & \quad i=1,\\
			\beta^{i-2}, & \quad 2\leq i\leq N.
		\end{array}
		\right.
	\end{align*}
	The infinity norm of the $kth$ power of $M$ is bounded as follows
	\begin{align*}
		\left\|M^{k}(\beta)\right\|_{\infty}\leq
		\left\{
		\begin{array}{ll}
			min
			\begin{Bmatrix}
				\left(\dfrac{1-|\beta|^{N-1}}{1-|\beta|}\right)^{k},
				\begin{pmatrix}
					N-1\\
					k
				\end{pmatrix}
			\end{Bmatrix}, &\quad |\beta|<1,\\
			|\beta|^{N-k-1}
			\begin{pmatrix}
				N-1\\
				k
			\end{pmatrix} ,&\quad |\beta| \geq 1.
		\end{array}
		\right.
	\end{align*}
\end{lemma}
\begin{lemma}\label{lemma6}\cite{GanderVandewalle2007}
	Let $\gamma,\eta\geq0$, a double indexed sequence \{$\delta_{n}^{k}$\} satisties $\delta_{n}^{k}\geq 0,\,\delta_{0}^{k}\geq 0$ and
	\begin{equation*}
		\delta_{n}^{k}\leq \gamma\delta_{n-1}^{k}+\eta\delta_{n-1}^{k-1},
	\end{equation*}
	for $n=0,1,\cdots,N$ and $k=0,1,\cdots,K$, then vector $\zeta^{k}=(\delta_{1}^{k},\delta_{2}^{k},\cdots,\delta_{N}^{k})^{T} $satisfies
	\begin{align*}
		\zeta^{k}\leq \eta M(\gamma)\zeta^{k-1}.
	\end{align*}
\end{lemma}
\noindent{\bf{Proof of Theorem 1.}} 
For all $n=0, \cdots, N$ and $k=0, \cdots, K$, denote the error $\varepsilon _{n}^{(k)}:=	\left\|u(t_{n})-u_{n}^{(k)}\right\|_{L_{2}(\Omega,\,\mathbb{H})}$. Since the exact solution  $u(t_{n})$  is chosen as the fine propagator $\mathcal{F}$, it can be written as
\begin{align}\label{eq16}
	u(t_{n})&=\mathcal{F}(t_{n-1},t_{n},u(t_{n-1}))\nonumber\\
	&=\mathcal{G} (t_{n-1},t_{n},u(t_{n-1}))+\mathcal{F} (t_{n-1},t_{n},u(t_{n-1}))-\mathcal{G} (t_{n-1},t_{n},u(t_{n-1})).
\end{align}
Subtracting the (\ref{eq16}) from (\ref{eq8}) and  using the notation of the residual operator (\ref{eq15}), we obtain
\begin{align*}
	\varepsilon _{n}^{(k)}	&\leq\left\|\mathcal{G} (t_{n-1},t_{n},u(t_{n-1}))-\mathcal{G} (t_{n-1},t_{n},u_{n-1}^{(k)})\right\|_{L_{2}(\Omega,\,\mathbb{H})}+\left\| \mathcal{R}(t_{n-1},t_{n},u(t_{n-1}))
	- \mathcal{R} (t_{n-1},t_{n},u_{n-1}^{(k-1)})\right\|_{L_{2}(\Omega,\,\mathbb{H})}\\
	&:=I_{1}+I_{2}.
\end{align*}
Firstly, we estimate $I_{1}$. Applying the stochastic exponential integrator (\ref{eq11}) for the coarse propagator $\mathcal{G}$, it holds that
\begin {align}
\mathcal{G} (t_{n-1},t_{n},u(t_{n-1}))&=\widehat S(\Delta T)u(t_{n-1})+\widehat S(\Delta T)F(u(t_{n-1}))\Delta T+\widehat S(\Delta T) B(t_{n-1})\Delta W_{n},\label{eq17}\\
\mathcal{G} (t_{n-1},t_{n},u_{n-1}^{(k)})&=\widehat S(\Delta T)u_{n-1}^{(k)}+\widehat S(\Delta T)F(u_{n-1}^{(k)})\Delta T+\widehat S(\Delta T)B(t_{n-1})\Delta W_{n}.\label{eq18}
\end{align}
Subtracting the above two formulas  leads to
\begin{align}\label{eq19}
I_{1}&=\left\|\widehat S(\Delta T)(u(t_{n-1})-u_{n-1}^{(k)})+\widehat S(\Delta T)(F(u(t_{n-1}))-F(u_{n-1}^{(k)}))\Delta T\right\|_{L_{2}(\Omega,\,\mathbb{H})}\nonumber \\
&\leq\left\|\widehat S(t)\right\|_{\mathcal{L(\mathbb{H})}}\left\|u(t_{n-1})-u_{n-1}^{(k)}\right\|_{L_{2}(\Omega,\,\mathbb{H})}+C\Delta T\left\|\widehat S(t)\right\|_{\mathcal{L(\mathbb{H})}}\left\|u(t_{n-1})-u_{n-1}^{(k)}\right\|_{L_{2}(\Omega,\,\mathbb{H})}\nonumber \\
&=(1+C\Delta T)\left\|u(t_{n-1})-u_{n-1}^{(k)}\right\|_{L_{2}(\Omega,\,\mathbb{H})}\nonumber \\
&=(1+C\Delta T)\varepsilon _{n-1}^{(k)},
\end{align}
which by the contraction property of semigroup and the global Lipschitz property of $F$.

Now it remains to estimate $I_{2}$. Applying exact solution integrator (\ref{eq12}) for fine progagator $\mathcal{F}$ leads to
\begin{align}
\mathcal{F}(t_{n-1},t_{n},u(t_{n-1}))&=\widehat S(\Delta T)u(t_{n-1})+\int_{0}^{\Delta T}\widehat S(\Delta T-s)F(U(t_{n-1},t_{n-1}+s,u(t_{n-1})))ds\nonumber\\
+\int_{0}^{\Delta T}\widehat S(\Delta T-s)B(t_{n-1})dW(s),\label{eq20}\\
\mathcal{F}(t_{n-1},t_{n},u_{n-1}^{(k-1)})&=\widehat S(\Delta T)u_{n-1}^{(k-1)}+\int_{0}^{\Delta T}\widehat S(\Delta T-s)F(V(t_{n-1},t_{n-1}+s,u_{n-1}^{(k-1)}))ds\nonumber\\
+\int_{0}^{\Delta T}\widehat S(\Delta T-s)B(t_{n-1})dW(s),\label{eq21}
\end{align}
where $U(t_{n-1},t_{n-1}+s,u)$ and $V(t_{n-1},t_{n-1}+s,u)$ denote the exact solution of system (\ref{eq4}) at time $t_{n-1}+s$ with the initial value $u$ and the initial time $t_{n-1}$.

Substituting the above equations and equations (\ref{eq17}) and (\ref{eq18})  into the residual operator (\ref{eq15}), we obtain
\begin{align*}
I_{2}&=\left\| \mathcal{R}(t_{n-1},t_{n},u(t_{n-1}))
- \mathcal{R} (t_{n-1},t_{n},u_{n-1}^{(k-1)})\right\|_{L_{2}(\Omega,\,\mathbb{H})}\\
&\leq \left\|\int_{0}^{\Delta T}\widehat S(\Delta T-s)[F(U(t_{n-1},t_{n-1}+s,u(t_{n-1})))-F(V(t_{n-1},t_{n-1}+s,u_{n-1}^{(k-1)}))]ds\right\|_{L_{2}(\Omega,\,\mathbb{H})}\\
&+\left\|\widehat S(\Delta T)[F(u(t_{n-1}))-F(u_{n-1}^{(k-1)})]\Delta T\right\|_{L_{2}(\Omega,\,\mathbb{H})}\\
&:=I_{3}+I_{4}.
\end{align*}
To get the estimation of $I_{3}$, by Lipschitz continuity property for $F$ and Lemma \ref{lemma4}, we derive
\begin{align}\label{eq22}
I_{3}&\leq  \int_{0}^{\Delta T} \left\|\widehat S(\Delta T-s)\right\|_{\mathcal{L(\mathbb{H})}}\left\|F(U(t_{n-1},t_{n-1}+s,u(t_{n-1})))-F(V(t_{n-1},t_{n-1}+s,u_{n-1}^{(k-1)}))\right\|_{L_{2}(\Omega,\,\mathbb{H})}ds\nonumber\\
&\leq C\int_{0}^{\Delta T} \left\|U(t_{n-1},t_{n-1}+s,u(t_{n-1}))-V(t_{n-1},t_{n-1}+s,u_{n-1}^{(k-1)})\right\|_{L_{2}(\Omega,\,\mathbb{H})}ds\nonumber\\
&\leq C\Delta T \left\|u(t_{n-1})-u_{n-1}^{(k-1)}\right\|_{L_{2}(\Omega,\,\mathbb{H})}.
\end{align}
As for $I_{4}$, using the contraction property of semigroup and Lipschitz continuity property for $F$ yields
\begin{align}\label{eq23}
I_{4}&\leq C\Delta T \left\|u(t_{n-1})-u_{n-1}^{(k-1)}\right\|_{L_{2}(\Omega,\,\mathbb{H})}.
\end{align}
From  (\ref{eq22}) and (\ref{eq23}), we know that
\begin{align}\label{eq24}
I_{2}&\leq C\Delta T\left\|u(t_{n-1})-u_{n-1}^{(k-1)}\right\|_{L_{2}(\Omega,\,\mathbb{H})}\nonumber\\
&=C\Delta T\varepsilon _{n-1}^{(k-1)}.
\end{align}
Combining ($\ref{eq19}$) and ($\ref{eq24}$) enables us to derive
\begin{align*}
\varepsilon _{n}^{(k)}\leq(1+C\Delta T)\varepsilon _{n-1}^{(k)}+C\Delta T\varepsilon _{n-1}^{(k-1)}.
\end{align*}
Let $\zeta^{k}=(\varepsilon_1^{k},\varepsilon_2^{k},\cdots,\varepsilon_{N}^{k})^{T}$ . It follows from Lemma \ref{lemma6} that 
\begin{align*}
\zeta^{k}&\leq C\Delta T M(1+C\Delta T)\zeta^{k-1}\\
&\leq C^{k} \Delta T^{k} M^{k}(1+C\Delta T)\zeta^{0}.
\end{align*}
Taking infinity norm and using Lemma \ref{lemma5} imply 
\begin{align*}
\sup\limits_{1\leq n \leq N}	\varepsilon _{n}^{(k)} &\leq(1+C\Delta T)^{N-k-1}C_{k}\Delta T^{k} C_{N-1}^{k}\sup\limits_{1\leq n \leq N}	\varepsilon _{n}^{(0)}\\
&\leq \frac{C_{k}}{k!} \Delta T^{k} \prod \limits_{j=1}^k (N-j)\sup\limits_{1\leq n \leq N}	\varepsilon _{n}^{(0)}.			 
\end{align*}
This completes the proof.\hfill $\Box$
\subsection{The stochastic exponential integrator as the fine propagator}
In this section, the error we considered  is the solution by the proposed algorithm and the reference solution generated by the fine propagator $\mathcal{F}$. To begin with, we define the reference solution as follows.
\begin{definition}
For all $n = 0,\ldots,N$, the reference solution is defined by  the fine propagator on each subinterval $[t_{n-1},t_{n}]$ 
\begin{align}\label{eq25}
	u_{n}^{ref}&=\mathcal{F}(t_{n-1},t_{n},u_{n-1}^{ref}),\\
	u_{0}^{ref}&=u_{0}.\nonumber
\end{align}
Precisely, 
\begin{align}\label{eq26}
	&u_{n-1,j}^{ref}=\mathcal{F}(t_{n-1,j-1},t_{n-1,j},u_{n-1,j-1}^{ref}), \, j=1,\cdots,J,\\
	&u_{n-1,0}^{ref}=u_{n-1}^{ref}.\nonumber
\end{align}
\end{definition}

\begin{theorem}\label{theorem2}
Let  \Cref{assump1}, \Cref{assump2} and \Cref{assump3} hold,  we apply the stochastic exponential integrator for coarse propagator $\mathcal{G}$ and the stochastic exponential integrator for fine propagator $\mathcal{F}$. Then we have the following convergence estimate for the fixed iteration number $k$ 
\end{theorem}
\begin{align}\label{eq27}
\sup\limits_{1\leq n \leq N}\left\|u_{n}^{(k)}-u_{n}^{ref}\right\|_{L_{2}(\Omega,\,\mathbb{H})}\leq\frac{C_{k}}{k!} \Delta T^{k} \prod \limits_{j=1}^k (N-j)\left\|u_{n}^{(0)}-u_{n}^{ref}\right\|_{L_{2}(\Omega,\,\mathbb{H})},
\end{align}
with a positive constant $C$ independent on $\Delta T$, where the parareal solution $u_{n}^{(k)}$ is defined in (\ref{eq8}) and the reference solution $u_{n}^{ref}$ is defined in (\ref{eq25}).

\noindent{\bf{Proof of Theorem 2.}}
For all $n = 0,\cdots,N$ and $k = 0,\cdots,K$, let the error be defined by 
$\varepsilon _{n}^{(k)}:=	\left\|u_{n}^{k}-u_{n}^{ref}\right\|_{L_{2}(\Omega,\,\mathbb{H})}$. 

Observe that the reference solution (\ref{eq25}) can be rewritten 
\begin{align}\label{eq28}
u_{n}^{ref}=\mathcal{F}(t_{n-1},t_{n},u_{n-1}^{ref})+\mathcal{G}(t_{n-1},t_{n},u_{n-1}^{ref})-\mathcal{G}(t_{n-1},t_{n},u_{n-1}^{ref}).
\end{align}
Combining the parareal algorithm form (\ref{eq8}) and the reference solution (\ref{eq28}) and using the notation of the residual operator (\ref{eq15}), the error can be written as
\begin{align*}
\varepsilon _{n}^{(k)}
&=\left\|\mathcal{G}(t_{n-1},t_{n},u_{n-1}^{(k)})-\mathcal{G}(t_{n-1},t_{n},u_{n-1}^{ref})+\mathcal{R}(t_{n-1},t_{n},u_{n-1}^{(k-1)})-\mathcal{R}(t_{n-1},t_{n},u_{n-1}^{ref})\right\|_{L_{2}(\Omega,\,\mathbb{H})}\\
&\leq\left\|\mathcal{G}(t_{n-1},t_{n},u_{n-1}^{(k)})-\mathcal{G} (t_{n-1},t_{n},u_{n-1}^{ref})\right\|_{L_{2}(\Omega,\,\mathbb{H})}+\left\| \mathcal{R}(t_{n-1},t_{n},u_{n-1}^{(k-1)})
- \mathcal{R}(t_{n-1},t_{n},u_{n-1}^{ref})\right\|_{L_{2}(\Omega,\,\mathbb{H})}\\
&:=I_{1}+I_{2}.
\end{align*}
Now we estimate $I_{1}$. Applying the stochastic exponential integrator (\ref{eq11}) for the coarse propagator $\mathcal{G}$, we obtain
\begin{align}\label{eq29}
\mathcal{G} (t_{n-1},t_{n},u_{n-1}^{ref})=\widehat S(\Delta T)u_{n-1}^{ref}+\widehat S(\Delta T)F(u_{n-1}^{ref})\Delta T+\widehat S(\Delta T)B(t_{n-1})\Delta W_{n}.
\end{align}
Subtracting the above formula (\ref{eq29}) from (\ref{eq18}), we have
\begin{align*}
I_{1}=\left\|\widehat S(\Delta T)(u_{n-1}^{(k)}-u_{n-1}^{ref})+\widehat S(\Delta T)[F(u_{n-1}^{(k)})-F(u_{n-1}^{ref})]\Delta T\right\|_{L_{2}(\Omega,\,\mathbb{H})}.
\end{align*}
Armed with contraction property of semigroup and Lipschitz continuity property of $F$ yield
\begin{align}\label{eq30}
I_{1}&\leq\left\|\widehat S(t)\right\|_{\mathcal{L(\mathbb{H})}}\left\|u_{n-1}^{(k)}-u_{n-1}^{ref}\right\|_{L_{2}(\Omega,\,\mathbb{H})}+
C\Delta T\left\|\widehat S(t)\right\|_{\mathcal{L(\mathbb{H})}} \left\|u_{n-1}^{(k)}-u_{n-1}^{ref}\right\|_{L_{2}(\Omega,\,\mathbb{H})}\nonumber\\
&\leq\left\|u_{n-1}^{(k)}-u_{n-1}^{ref}\right\|_{L_{2}(\Omega,\,\mathbb{H})}+
\Delta T C\left\|u_{n-1}^{(k)}-u_{n-1}^{ref}\right\|_{L_{2}(\Omega,\,\mathbb{H})}\nonumber\\
&\leq (1+C\Delta T )\left\|u_{n-1}^{(k)}-u_{n-1}^{ref}\right\|_{L_{2}(\Omega,\,\mathbb{H})}.
\end{align}
As for $I_{2}$, regarding the estimation of the residual operator, we need to resort to its directional derivatives. Due to formula ($\ref{eq15}$), the derivatives can be given by
\begin{align}\label{eq31}
D\mathcal{R}(t_{n-1},t_{n},u).h:=D\mathcal{F}(t_{n-1},t_{n},u).h-D\mathcal{G}(t_{n-1},t_{n},u).h.
\end{align}
One the one hand, since the stochastic exponential scheme is chosen as the fine propagator (\ref{eq13}) with time step-size $\Delta t$, we obtain
\begin{align*}
\left\{
\begin{array}{ll}
	u _{n,j+1}=\widehat S(\Delta t)u _{n,j} +\Delta t \widehat S(\Delta t)F(u_{n,j})+\widehat S(\Delta t)B(t_{n,j})\Delta W_{n,j},\quad j\in1,\cdots,J-1,\\
	u _{n,0}=u .
\end{array}
\right.
\end{align*}
Denote  $D(u_{n,j}).h:=\eta^{h}_{n,j} $ for $j \in 0,\cdots,J$. Then  taking the direction derivatives for above equation yields
\begin{align*}
\left\{
\begin{array}{ll}
	\eta _{n,j+1}^{h}=\widehat S(\Delta t)\eta _{n,j}^{h} +\Delta t \widehat S(\Delta t)DF(u _{n,j}).\eta _{n,j}^{h},\quad j \in1,\cdots,J-1,\\
	\eta _{n,0}^{h}=h .
\end{array}
\right.
\end{align*}
Based on the the form of semigroup  $\{\widehat S(t)=e^{ t(M-\sigma Id)},t\geq 0\}$, we have the following recursion formula
\begin{align*}
\eta _{n,J}^{h}&=e^{J\Delta t(M-\sigma Id)}\eta _{n,0}^{h}+\Delta t\sum\limits_{j=0}\limits^{J-1}e^{(J-j)\Delta t(M-\sigma Id)}DF(u _{n,j}).\eta _{n,j}^{h}\\
&=e^{\Delta T(M-\sigma Id)}h+\Delta t\sum\limits_{j=0}\limits^{J-1}e^{(J-j)\Delta t(M-\sigma Id)}DF(u _{n,j}). \eta _{n,j}^{h}.
\end{align*}
Applying the discrete Gronwall lemma yields the following inequality
\begin{align}\label{eq32}
\sup\limits_{1\leq n \leq N}\left\|\eta _{n,j}^{h}\right\|_{L_{2}(\Omega,\,\mathbb{H})}\leq C\left\|h\right\|_{L_{2}(\Omega,\,\mathbb{H})}.
\end{align}
Moreover, the derivative of $\mathcal{F}(t_{n-1},t_{n},u)$ can be writen by $D\mathcal{F}(t_{n-1},t_{n},u).h=D(u_{n,J}).h=\eta_{n,J}^{h}$, where $J\Delta t=\Delta T$, that is, one gets
\begin{align}\label{eq33}
D\mathcal{F}(t_{n-1},t_{n},u).h=e^{\Delta T(M-\sigma Id)}h+\Delta t\sum\limits_{j=0}\limits^{J-1}e^{(J-j)\Delta t(M-\sigma Id)}DF(u_{n,j}). \eta _{n,j}^{h}.
\end{align}
On the other hand, since the stochastic exponential scheme is chosen as the coarse propagator $\mathcal{G}$, taking the direction derivative for $u$ of formula (\ref{eq11}) leads to
\begin{align}\label{eq34}
D\mathcal{G}(t_{n-1},t_{n},u).h=e^{\Delta T(M-\sigma Id)}h+\Delta Te^{\Delta T(M-\sigma Id)}DF(u).h.
\end{align}
Substituting formula (\ref{eq33}) and (\ref{eq34}) into formula (\ref{eq31}), we obtain
\begin{align*}
\left\|D\mathcal{R}(t_{n-1},t_{n},u).h\right\|_{L_{2}(\Omega,\,\mathbb{H})}&= \left\|D\mathcal{F}(t_{n-1},t_{n},u).h-D\mathcal{G}(t_{n-1},t_{n},u).h\right\|_{L_{2}(\Omega,\,\mathbb{H})}\\
&=\left\|\Delta t\sum\limits_{j=0}\limits^{J-1}e^{(J-j)\Delta t(M-\sigma Id)}DF(u _{n,j}). \eta _{n,j}^{h}+\Delta Te^{\Delta T(M-\sigma Id)}DF(u).h \right\|_{L_{2}(\Omega,\,\mathbb{H})}\\
&\leq\left\|\Delta t\sum\limits_{j=0}\limits^{J-1}e^{(J-j)\Delta t(M-\sigma Id)}DF(u_{n,j}).\eta _{n,j}^{h} \right\|_{L_{2}(\Omega,\,\mathbb{H})}+\left\|\Delta Te^{\Delta T(M-\sigma Id)}DF(u).h\right\|_{L_{2}(\Omega,\,\mathbb{H})}.
\end{align*}
Utilizing the bounded derivatives condition of $F$, we get
\begin{align*}
\left\|D\mathcal{R}(t_{n-1},t_{n},u).h\right\|_{L_{2}(\Omega,\,\mathbb{H})}
\leq\Delta t\sum\limits_{j=0}\limits^{J-1}\left\|e^{(J-j)\Delta t(M-\sigma Id)}\right\|_{\mathcal{L(\mathbb{H})}}\left\|\eta _{n,j}^{h}\right\|_{L_{2}(\Omega,\,\mathbb{H})}+C\Delta T\left\|e^{\Delta T(M-\sigma Id)}\right\|_{\mathcal{L(\mathbb{H})}}\left\|h\right\|_{L_{2}(\Omega,\,\mathbb{H})}.
\end{align*}
Using the contraction property of semigroup, we have
\begin{align*}
\sup\limits_{1\leq n \leq N}\left\|D\mathcal{R}(t_{n-1},t_{n},u).h\right\|_{L_{2}(\Omega,\,\mathbb{H})}
\leq\Delta t\sum\limits_{j=0}\limits^{J-1}\sup\limits_{1\leq n \leq N}\left\|\eta _{n,j}^{h}\right\|_{L_{2}(\Omega,\,\mathbb{H})}+C\Delta T\left\|h\right\|_{L_{2}(\Omega,\,\mathbb{H})}.
\end{align*}
Substituting the Gronwall inequality ($\ref{eq32}$) into the above inequality leads to
\begin{align*}
\sup\limits_{1\leq n \leq N}\left\|D\mathcal{R}(t_{n-1},t_{n},u).h\right\|_{L_{2}(\Omega,\,\mathbb{H})}
&\leq\Delta t J C\left\|h\right\|_{L_{2}(\Omega,\,\mathbb{H})}+C\Delta T\left\|h\right\|_{L_{2}(\Omega,\,\mathbb{H})}\\
&\leq C\Delta T\left\|h\right\|_{L_{2}(\Omega,\,\mathbb{H})}.
\end{align*}
In conclusion, it holds that
\begin{align}\label{eq35}
\sup\limits_{1\leq n \leq N}\left\|\mathcal{R}(t_{n-1},t_{n},u_{2})-\mathcal{R}(t_{n-1},t_{n},u_{1})\right\|_{L_{2}(\Omega,\,\mathbb{H})}\leq C \Delta T\left\|u_{2}-u_{1}\right\|_{L_{2}(\Omega,\,\mathbb{H})},\,\forall u_{1}, u_{2}\in \mathbb{H}.
\end{align}
Substituting $u_{n-1}^{(k-1)}$ and $u_{n-1}^{ref}$ into above formula derives lipschitz continuity property of the residual operator
\begin{align}\label{eq36}
I_{2}&=\left\|\mathcal{R}(t_{n-1},t_{n},u_{n-1}^{(k-1)})-\mathcal{R}(t_{n-1},t_{n},u_{n-1}^{ref})\right\|_{L_{2}(\Omega,\,\mathbb{H})}\leq C\Delta T \left\|u_{n-1}^{(k-1)}-u_{n-1}^{ref}\right\|_{L_{2}(\Omega,\,\mathbb{H})}.
\end{align}
Combining (\ref{eq30}) and (\ref{eq36}), we have
\begin{align*}
\varepsilon _{n}^{(k)}
&\leq (1+C\Delta T )\left\|u_{n-1}^{(k)}-u_{n-1}^{ref}\right\|_{L_{2}(\Omega,\,\mathbb{H})}+C\Delta T \left\|u_{n-1}^{(k-1)}-u_{n-1}^{ref}\right\|_{L_{2}(\Omega,\,\mathbb{H})}\\
&=(1+C\Delta T)\varepsilon _{n-1}^{(k)}+C\Delta T \varepsilon _{n-1}^{(k-1)}.
\end{align*}
According to Lemma \ref{lemma5} and Lemma \ref{lemma6}, it yields to
\begin{align*}
\sup\limits_{1\leq n \leq N}	\varepsilon _{n}^{(k)} &\leq(1+C\Delta T)^{N-k-1}C_{k}\Delta T^{k} C_{N-1}^{k}\sup\limits_{1\leq n \leq N}	\varepsilon _{n}^{(0)}\\
&\leq \frac{C_{k}}{k!} \Delta T^{k} \prod \limits_{j=1}^k (N-j)\sup\limits_{1\leq n \leq N}	\varepsilon _{n}^{(0)},		 
\end{align*}
which leads to the final result
\begin{align*}
\sup\limits_{1\leq n \leq N}\left\|u_{n}^{(k)}-u_{n}^{ref}\right\|_{L_{2}(\Omega,\,\mathbb{H})}\leq\frac{C_{k}}{k!} \Delta T^{k} \prod \limits_{j=1}^k (N-j)\left\|u_{n}^{(0)}-u_{n}^{ref}\right\|_{L_{2}(\Omega,\,\mathbb{H})}.
\end{align*}\hfill $\Box$
\begin{remark}
We can summarise Lipschitz continuity property of the residual operator $\mathcal{R}(t_{n-1},t_{n},u)$: there exists $C\in(0,\infty)$ such that for $\Delta T\in\left(0,1\right]$ and  $u_{1},u_{2}\in \mathbb{H}$, we have
\begin{align*}
	\sup\limits_{1\leq n \leq N}\left\|\mathcal{R}(t_{n-1},t_{n},u_{2})-\mathcal{R}(t_{n-1},t_{n},u_{1})\right\|_{L_{2}(\Omega,\,\mathbb{H})}\leq C \Delta T\left\|u_{2}-u_{1}\right\|_{L_{2}(\Omega,\,\mathbb{H})}.
\end{align*}
\end{remark} 
\begin{remark}
When we fix the iteration number $k$, the convergence rate will be $\mathcal{O}(\Delta T)^{k}$.
\end{remark}
\begin{remark}
The error between the reference solution $u_{n}^{ref}$ by the fine propagator defined in (\ref{eq25}) and the exact solution $u(t_{n})$ defined in (\ref{eq9}) do not affect the convergence rate of the parareal algorithm, due to
\begin{align*}
	\sup\limits_{1\leq n \leq N}\left\|u_{n}^{ref}-u(t_{n})\right\|_{L_{2}(\Omega,\,\mathbb{H})}\leq C(1+\left\|u_0\right\|_{L_{2}(\Omega,\,\mathbb{H})})\Delta t,
\end{align*}
Therefore, it is sufficient to study the convergence order of  the error between $u_{n}^{(k)}$ and $u_{n}^{ref}$.
\begin{proposition}\cite{Cohenetal2020}
	{\bf(Uniform boundedness of reference solution $u_{n}^{ref}$)}. There exists a constant $C\in (0,\infty)$ such that 
	\begin{equation*}
		\sup\limits_{t\in[0,T]}\left\|u_{n}^{ref}\right\|_{L_{2}(\Omega,\,\mathbb{H})}\leq C(1+\left\|u_0\right\|_{L_{2}(\Omega,\,\mathbb{H})}).	
	\end{equation*}
\end{proposition}
\begin{proposition}
	{\bf(Uniform boundedness of  parareal algorithm solution $u_{n}^{(k)}$)}. There exists a constant $C\in (0,\infty)$ such that
	\begin{equation*}
		\sup\limits_{t\in[0,T]}\left\|u_{n}^{(k)}\right\|_{L_{2}(\Omega,\,\mathbb{H})}\leq C(1+\left\|u_0\right\|_{L_{2}(\Omega,\,\mathbb{H})}).	
	\end{equation*}
\end{proposition}

\end{remark}
\section{Numerical  experiments}

This section is devoted to investigating the convergence result with several parameters and the effect of the scale of noise on numerical solutions. Since the parareal algorithm in principle is a temporal algorithm, and the spatial discretization is not our focus in this article, we perform finite difference method to discretize spatially. 

The mean-square error is used as
\begin{align*}
L_2=(\varepsilon\|\vec{E}(t_n)-\vec{E}^{(k)}_{n}\|^2+\mu\|\vec{H}(t_n)-\vec{H}^{(k)}_{n}\|^2)^{\frac{1}{2}}.
\end{align*}
\subsection{Convergence}
\subsubsection{One-dimensional transverse magnetic wave}
We first consider the stochastic Maxwell equations with 1-D transverse magnetic wave driven by the standard Brownian motion
\begin{align*}
\left\{\begin{array}{l}
	\dfrac{\partial E_{z}}{\partial t}=\dfrac{1}{\epsilon}\dfrac{\partial H_{y}}{\partial x}-\sigma E_{z}+\lambda_{1}\cdot \dot{W}, \\[4mm]
	\dfrac{\partial H_{y}}{\partial t}=\dfrac{1}{\mu}\dfrac{\partial E_{z}}{\partial x}-\sigma H_{y}+\lambda_{2}\cdot \dot{W},
\end{array}\right.
\end{align*}
by providing initial conditions
\begin{align*}
E_{z}(x, 0)=\sin (x), \quad H_{y}(x,0)=-\sqrt{\frac{\epsilon}{\mu}} \sin (x),
\end{align*}
for $t\in [0,1]$, $x\in [0,2\pi]$ and $\dot{W} = \dot{B}$.
\begin{figure}[htbp]
	\centerline{\includegraphics[width=4in]{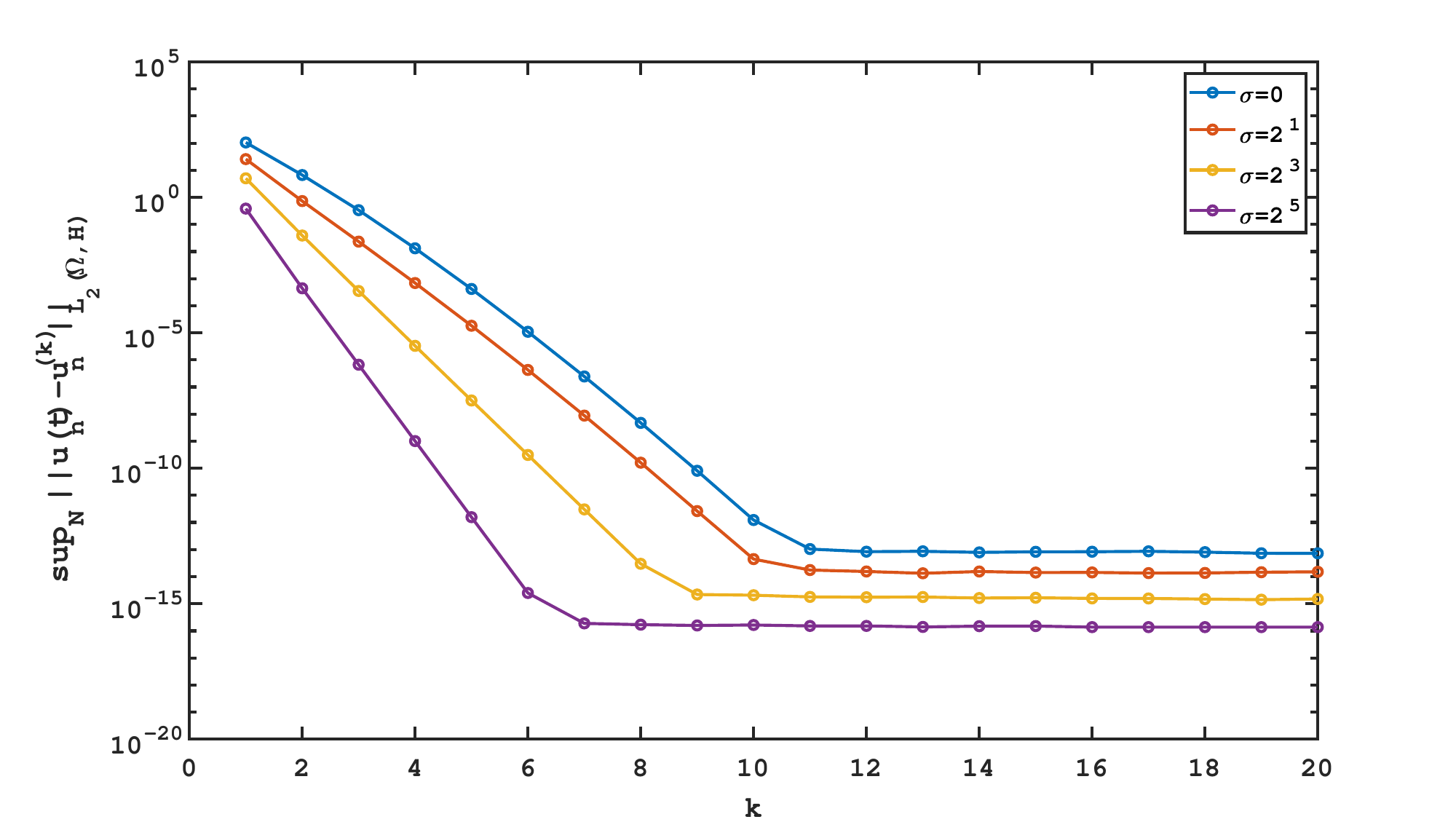}}
	\vspace*{8pt}
	\caption{Convergence of 1D case vs. interation number $k$ for different values of $\sigma=0,2^1,2^3,2^5$.}
	\label{fig1}
\end{figure}

The parameters are normalized to $\varepsilon =1$, $\mu =1$ and $\lambda_{1}=\lambda_{2}=1$. We apply the parareal algorithm to solve the numerical solution with the fine step-size  $\Delta t= 2^{-8}$ and the coarse step-size $\Delta T= 2^{-6}$. The spatial mesh grid-size $\Delta x= \Delta y=2\pi/100$. Figure \ref{fig1} demonstrates the evolution of the mean-square error $(\sup_{1\leq n \leq N}E\|u_{n}^{(k)}-u_{n}^{ref}\|^2)^{\frac{1}{2}}$ with the iteration number $k$. From the Figure.\ref{fig1}, we observe that the damping term speeds up the convergence of the numerical solutions and the error approaches $10^{-12}$ after $k = 11$ at least nearly, which shows that the proposed algorithm converges. 

\begin{remark}
From a numerical analysis point of view, the inclusion of damping coefficients usually accelerates the convergence of numerical solutions by suppressing oscillations and instability, resulting in a faster steady state or desired precision. However, too small damping may not be enough to accelerate the convergence rate and may even introduce instability.
\end{remark}
\begin{figure}[htbp]
	\centerline{\includegraphics[width=4in,height=3in]{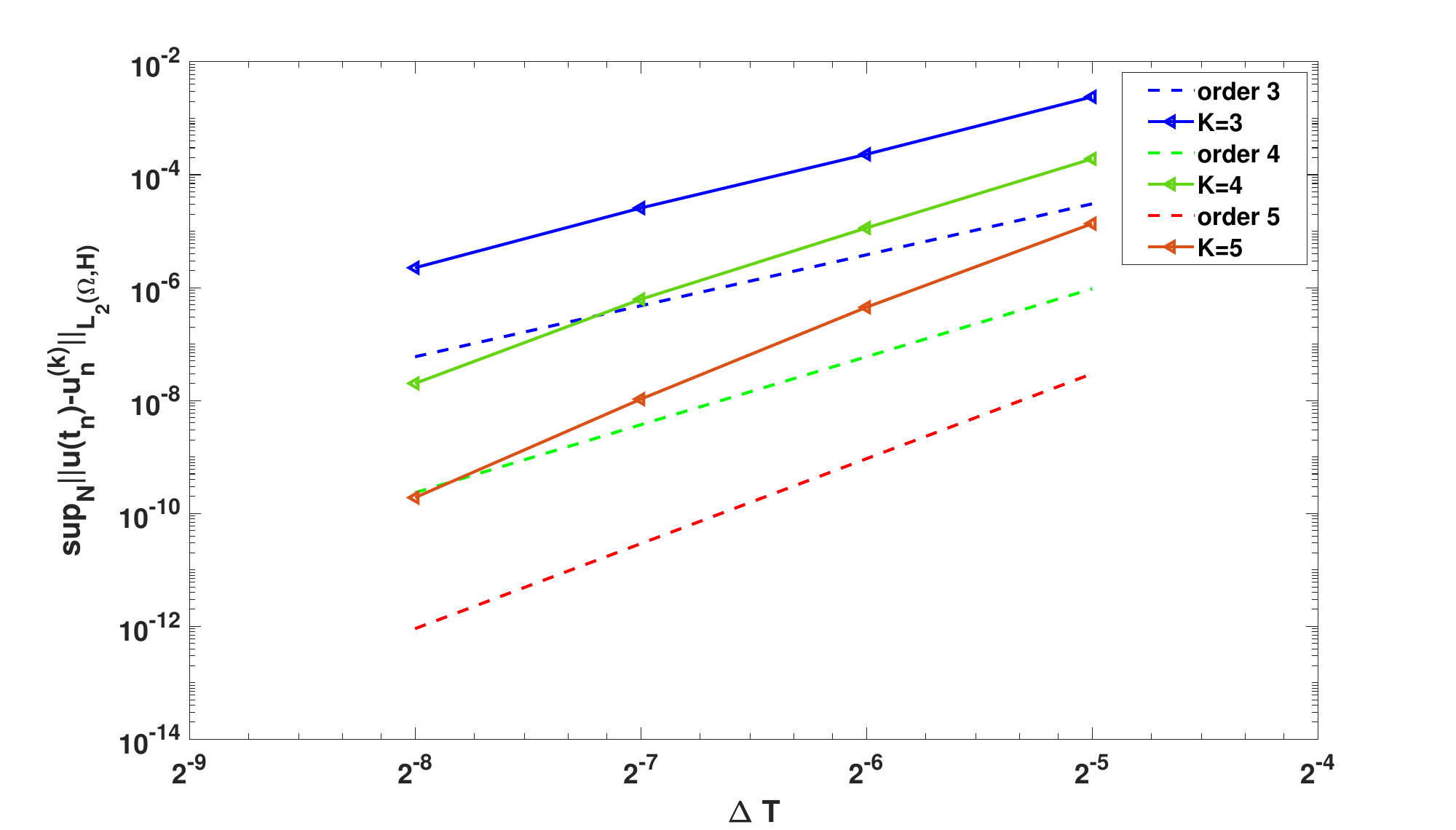}}
	\vspace*{8pt}
	\caption{Mean-square order of 1D case with respect to $\Delta T = 2^{-i}, i = 5,6,7,8.$}
	\label{fig2}
\end{figure}
Subsequently, we choose the damping coefficient $\sigma=2$ to calculate the convergence order of the proposed  parareal algorithm. We compute the numerical solution with the fine step-size  $\Delta t= 2^{-10}$ and the coarse step-size $\Delta T=2^{-5},2^{-6},2^{-7},2^{-8}$. Figure \ref{fig2} reports the  convergence order of the parareal algorithm with the iteration number $k = 3, 4, 5$. It is clearly shown that the mean-square convergence order always increases as the iteration number $k$ increases.

\subsubsection{Two-dimensional  transverse magnetic  waves}
We consider the stochastic Maxwell equations  with 2-D transverse magnetic polarization driven by trace class noise
\begin{align}\label{eq37}
\left\{
\begin{array}{ll}
	\dfrac{\partial E_{z}} {\partial t}=\dfrac{1} {\varepsilon}(\dfrac{\partial H_{y}} {\partial x}-\dfrac{\partial H_{x}} {\partial y})-\sigma E_{z}+\lambda_{1}\cdot \dot{W},\\[4mm]
	\dfrac{\partial H_{x}}{\partial t}=-\dfrac{1}{\mu}\dfrac{\partial E_{z}}{\partial y}-\sigma H_{x}+\lambda_{2}\cdot \dot{W},\\ [4mm]
	\dfrac{\partial H_{y}}{\partial t}=\dfrac{1}{\mu}\dfrac{\partial E_{z}}{\partial x} -\sigma H_{y}+\lambda_{2}\cdot \dot{W},
\end{array}\right.
\end{align}
by providing initial conditions
\begin{align*}
E_{z}(x,y,0)&=\sin(3\pi x)\sin(4\pi y) ,\\[2mm]
H_{x}(x,y,0)&=-0.8\cos(3\pi x)\sin(4\pi y) ,\\[2mm] 
H_{y}(x,y,0)&=-0.6\sin(3\pi x)\sin(4\pi y),
\end{align*}
for $t\in [0,1]$ and $(x,y)\in D=[0,2\pi]\times[0,2\pi]$. Following the formula (\ref{eq2}), we choose $e_{n}(x)=\sqrt{2/a}\sin\left(n\pi x/a\right)$ and $\lambda_{n}=n^{-(2r+1+\delta )}$, for some $\delta >0$ and $r\geq 0$. In this case, $Tr\left(Q\right)=\sum_{n=1}^{\infty}\lambda_{n}<\infty$. We construct the Wiener process as follows \citep{Chenetal2023}
\begin{align*}
W(t)=\sqrt{\dfrac{2}{a}}\sum_{n=1}^{\infty}n^{-(\frac{2r+1+\delta }{2})}\sin\left(\frac{n\pi x}{a}\right)\beta_{n}(t),		
\end{align*}
with $a=2,\,r=0.5$ and $\delta =0.001$.
\begin{figure}[htbp]
	\centerline{\includegraphics[width=4in]{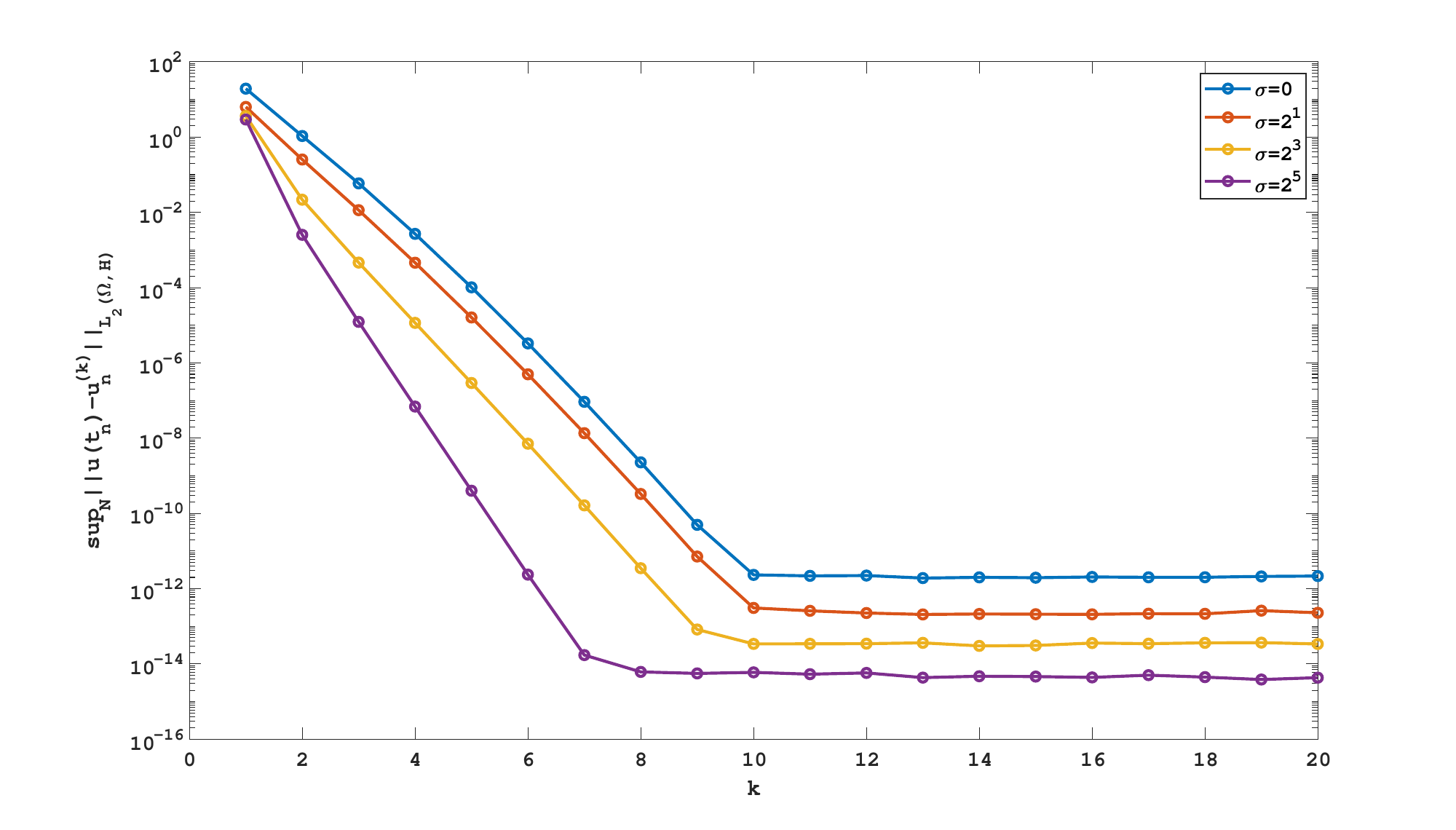}}
	\vspace*{8pt}
	\caption{Convergence of  2D case with interation number $k$ for different values of $\sigma=0,2^1,2^3,2^5$.}
	\label{fig3}
\end{figure}

Firstly, the parameters are normalized to $\varepsilon =1$, $\mu =1$ and $\lambda_{1}=\lambda_{2}=1$.
We take  the fine step-size $\Delta t= 2^{-8}$, the coarse step-size $\Delta T= 2^{-6} $ and the spatial mesh grid-size $\Delta x= \Delta y=2\pi/100$. Figure  \ref{fig3} demonstrates the evolution of the mean-square error $(\sup_{1\leq n \leq N}E\|u_{n}^{(k)}-u_{n}^{ref}\|^2)^{\frac{1}{2}}$ with iteration number $k$. From the Figure \ref{fig3}, we observe that the error approaches $10^{-11}$ after $k = 10$ nearly, which shows that the proposed algorithm converges.

\begin{remark}
	In numerical simulation, the introduction of damping terms and the selection of parameters need to be careful to ensure the accuracy and physical authenticity of simulation results. Excessive damping may lead to excessive attenuation, thus affecting the accuracy of simulation results.
\end{remark}

\begin{figure}[htbp]
	\centerline{\includegraphics[width=4in,height=3in]{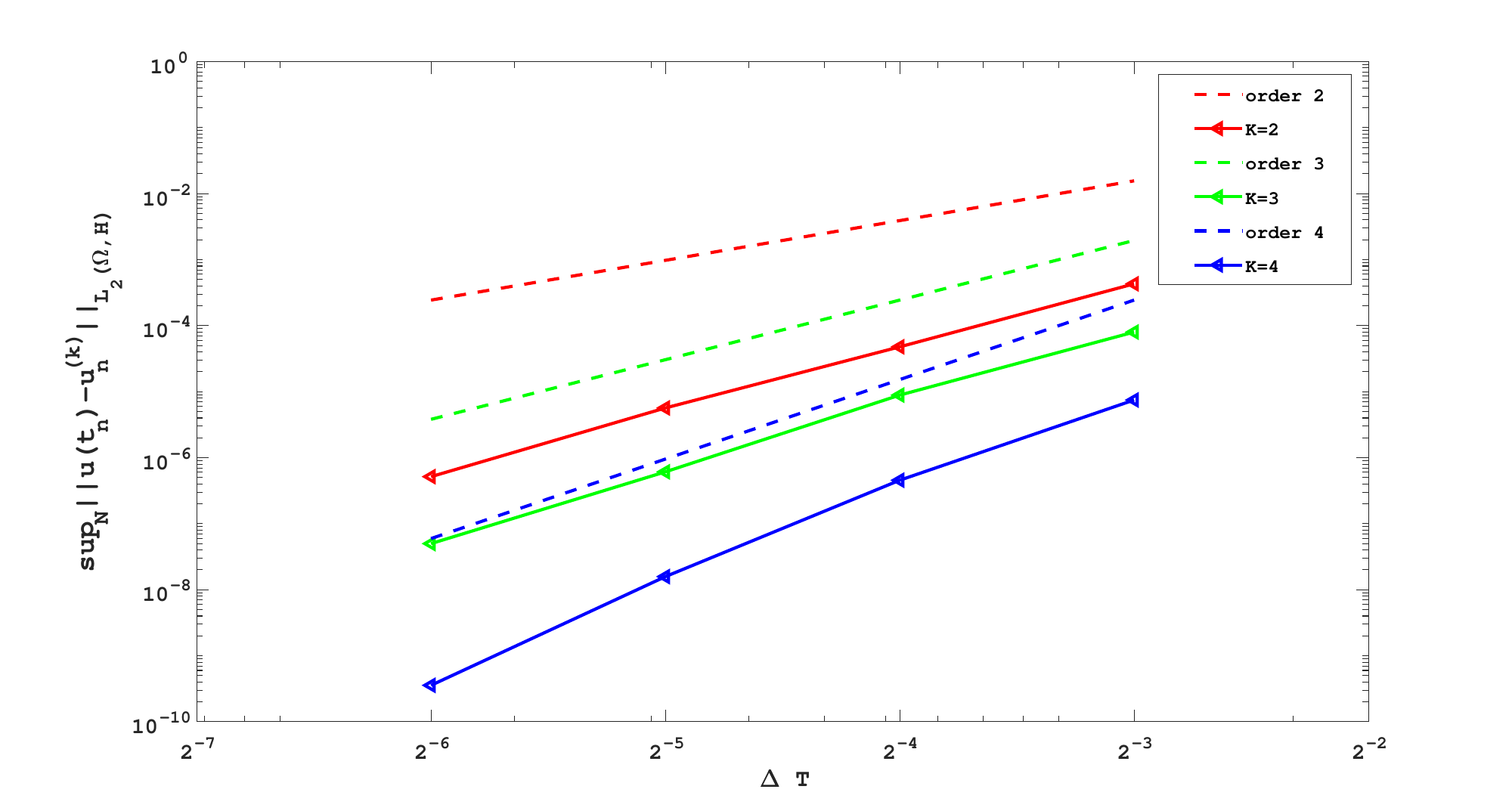}}
	\vspace*{8pt}
	\caption{Mean-square order of 2D case  with respect to $\Delta T = 2^{-i}, i = 3,4,5,6.$}
	\label{fig4}
\end{figure}

Secondly, in order to investigate the relationship between the convergence order and the iteration number, we choose the damping coefficient $\sigma=2$ to calculate the convergence order of the proposed  algorithm as taking the different iteration number $k$. We compute the numerical solution with the fine step-size  $\Delta t= 2^{-8}$ and the coarse step-size $\Delta T=2^{-3},2^{-4},2^{-5},2^{-6}$. Figure \ref{fig4} reports the  convergence order of the proposed algoritnumerical errorhm with the iteration number $k = 2, 3, 4$. Indeed, the numerical experiments reveal that the convergence order of the proposed algorithm increases as the iteration number $k$ increases.

\subsection{Impact of the scale of noise}

We consider the stochastic Maxwell equations  with 2-D transverse magnetic polarization (\ref{eq37}). The parameters are normalized to $\varepsilon =1$, $\mu =1$ and we take the fine step-size $\Delta t= 2^{-5}$, the coarse step-size $\Delta T= 2^{-3} $ and the spatial mesh grid-size $\Delta x= \Delta y=1/50$. In order to show 
the impact of the scale of noise on the numerical solution, we perform numerical simulations with four scales of noise $\lambda_{1}=\lambda_{2}=0, 2^1, 2^3, 2^5$ and choose the damping coefficient $\sigma=2^3$.
\begin{figure}[ht]
\centerline{\includegraphics[width=7in]{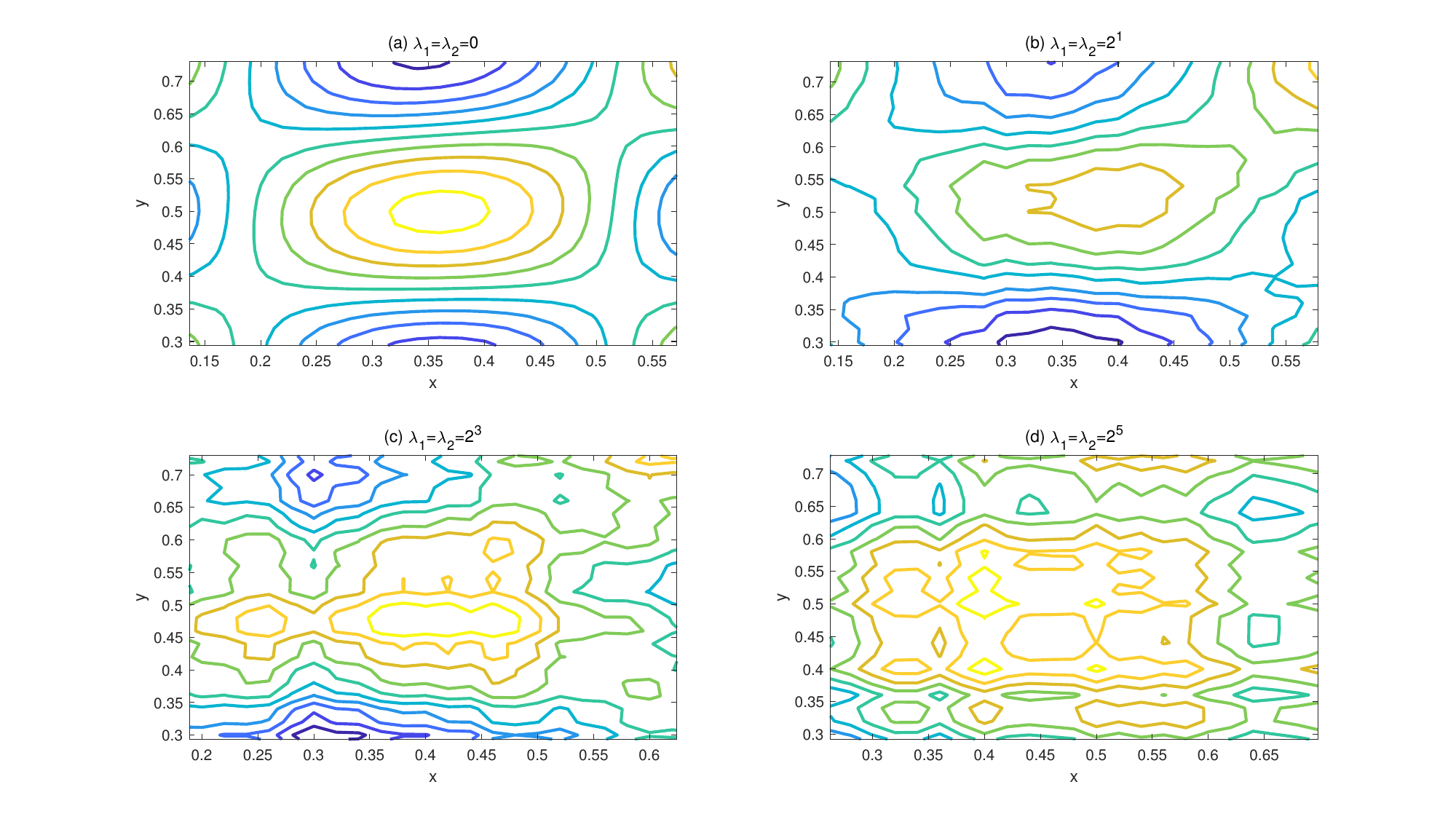}}
\vspace*{8pt}
\caption{10 Contour of $E_{z}(x,y)$ with different sizes of noise $\lambda_{1}=\lambda_{2}=0, 2^1, 2^3, 2^5$ in the time $T=1$.}
\label{fig5}
\end{figure}

\begin{figure}[ht]
\centerline{\includegraphics[width=7in]{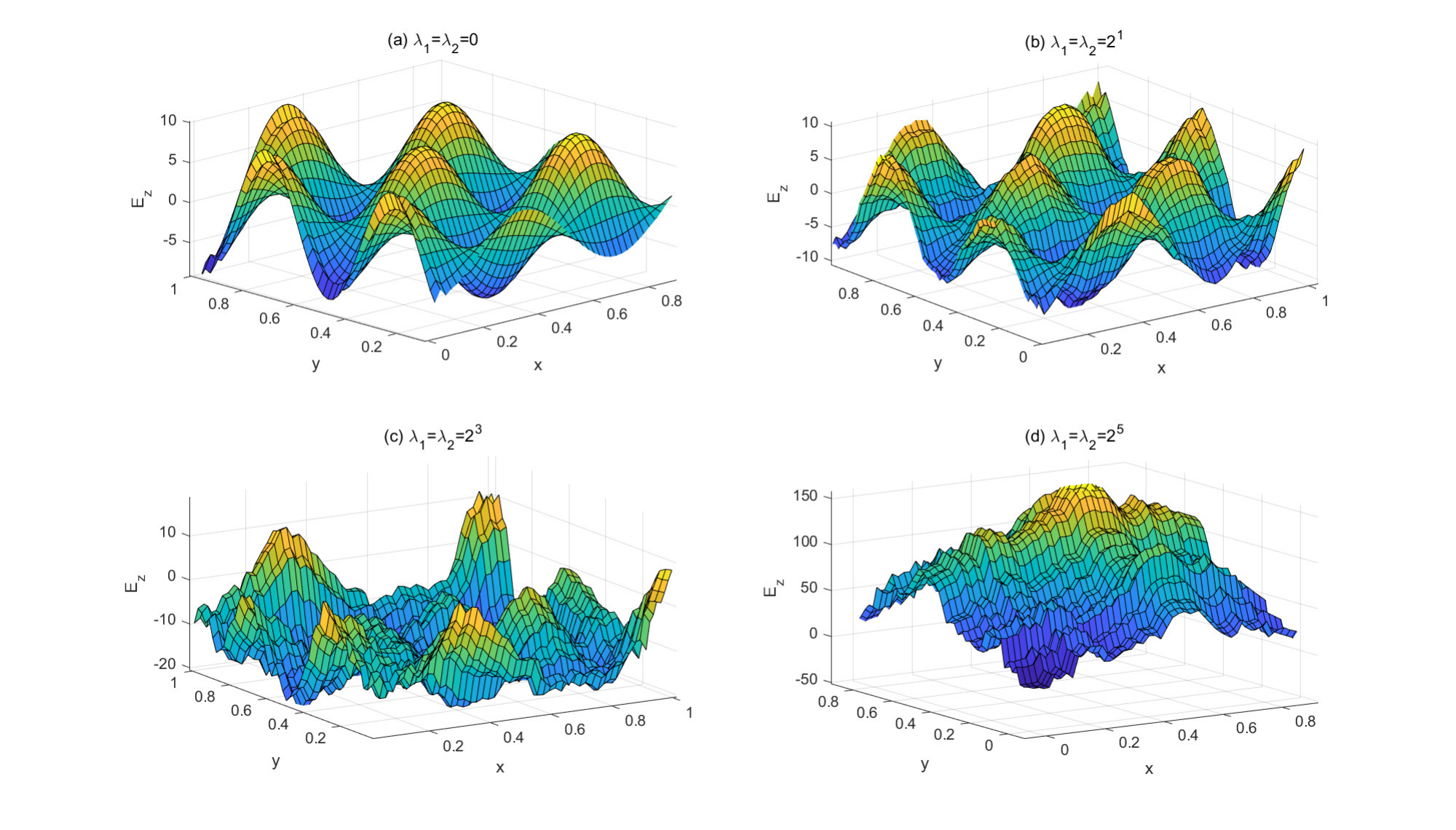}}
\vspace*{8pt}
\caption{$E_{z}(x,y)$ with different sizes of noise $\lambda_{1}=\lambda_{2}=0, 2^1, 2^3, 2^5$ in the time $T=1$.}
\label{fig6}
\end{figure}

Figure \ref{fig5} shows the 10 Contour plots of the numerical solution $E_{z}(x,y)$  with different scales of noise and Figure \ref{fig6} shows the electric field wave forms $E_{z}(x,y)$ with different scales of noise. Comparing with deterministic case (a) of Figure \ref{fig5} and Fig.\ref{fig6}, we can find that the oscillator of the wave forms (b-d) of Figure \ref{fig5} and Figure \ref{fig6} becomes more and more violent as the scale of the noise increases, i.e., from (a-d) of Figure \ref{fig5} and Figure \ref{fig6} it can be observed  that the perturbation of the numerical solutions becomes more and more apparent as the scale of the noise increases.
\section{Conclusion}	
In this paper, we study the strong convergence analysis of the parareal algorithms for stochastic Maxwell equations with damping term driven by additive noise. Firstly the stochastic exponential scheme is chosen as the coarse propagator and the exact solution scheme is chosen as the fine propagator. And we propose our numerical schemes and establish the mean-square convergence estimate. Secondly, both  the coarse propagator and the fine propagator choose the stochastic exponential scheme. Meanwhile, the error we considered in this section is the distance between the solution computed by the parareal algorithm and the reference solution generated by the fine propagator. It is shown that the convergence order of the proposed algorithms is linearly related to the iteration number $k$.
At last, One- and two-dimensional numerical examples are performed to demonstrate convergence analysis with respect to damping coefficient and noise scale. One key idea from the proofs of two convergence results is that the residual operator in Theorem \ref*{theorem2} is related to Lipschitz continuity properties, whereas Theorem \ref*{theorem1} concerns the integrability of the exact solution. The future works will include the study for the parareal algorithms for the stochastic Maxwell equations driven by multiplicative noise and other choices of integrators as the coarse and fine propagators.
\section*{Acknowledgments}
The authors would like to express their appreciation to the referees for their useful comments and the editors. Liying Zhang is supported by 
the National Natural Science Foundation of China (No.\,11601514 and No.\,11971458), the Fundamental Research Funds for the Central Universities (No.\,2023ZKPYL02 and No.\,2023JCCXLX01) and the Yueqi Youth Scholar Research Funds for the China University of Mining and Technology-Beijing (No.\,2020YQLX03).
\bibliography{refs}
\bibliographystyle{plain}
\end{document}